\documentclass[letterpaper, double spaced,10pt]{article}
\title{Riemann flow and Riemann wave via \\bialternate product Riemannian metric}
\author{Constantin Udri\c ste\footnote{Prof. Dr., University  Politehnica of Bucharest, Faculty of Applied Sciences,
Department of Mathematics-Informatics, Splaiul Independen\c tei 313, Bucharest, 060042, Romania.\newline
E-mail address: udriste@mathem.pub.ro, anet.udri@yahoo.com}}

\begin{document}
\maketitle

\newcommand{\area}{\hbox{area}}
\newcommand{\ome}{\Omega}

\newcommand{\omm}{\Omega}
\newcommand{\grad}{\hbox{$\:$grad \,}}
\newcommand{\divv}{\hbox{$\:$div \,}}
\newcommand{\rot}{\hbox{$\:$rot \,}}
\newcommand{\im}{\hbox{$\:$Im }}
\newcommand{\re}{\hbox{$\:$Re }}
\newcommand{\sh}{\hbox{$\:$sh\, }}
\newcommand{\ch}{\hbox{$\:$ch\,}}
\newcommand{\om}{\omega}
\newcommand{\br}{\hbox{\bf R}}
\newcommand{\bc}{\hbox{\bf C}}
\newcommand{\bz}{\hbox{\bf Z}}
\newcommand{\pp}{\prime}
\newcommand{\ty}{\infty}
\newcommand{\di}{\displaystyle}
\newcommand{\va}{\varphi}
\newcommand{\si}{\sigma}
\newcommand{\ga}{\gamma}
\newcommand{\gaa}{\Gamma}
\newcommand{\na}{\nabla}
\newcommand{\te}{\theta}
\newcommand{\ld}{\ldots}
\newcommand{\ov}{\over}
\newcommand{\ri}{\Rightarrow}
\newcommand{\noa}{\noalign{\medskip}}
\newcommand{\la}{\lambda}
\newcommand{\su}{\subset}
\newcommand{\qu}{\quad}
\newcommand{\fo}{\forall}
\newcommand{\al}{\alpha}
\newcommand{\lll}{\Leftrightarrow}
\newcommand{\be}{\beta}
\newcommand{\ep}{\varepsilon}
\newcommand{\pa}{\partial}
\newcommand{\ti}{\times}
\newcommand{\dd}{\Delta}
\newcommand{\sgg}{\Sigma}
\newcommand{\de}{\delta}

\newcommand\C{{\,I\!\!\!\!\!\: C}}
\newcommand{\ds}{\displaystyle}
\newcommand{\ol}{\overline}
\newcommand{\w}{\widetilde}
\newcommand{\ul}{\underline}
\newcommand{\s}{\stackrel}
\newcommand{\p}{\prime}
\def\d{\displaystyle\frac}
\def\pa{\partial}
\def\wi{\widetilde}
\def\ov{\over}
\def\ri{\rightarrow}
\def\ba{\begin{array}}
\def\ea{\end{array}}
\def\Ker{{\rm{Ker}}\,}
\def\var{\varepsilon}
\def\fii{\varphi}
\def\Lam{\Lambda}
\def\lam{\lambda}
\def\ii{\^\i}
\def\dx{\dot x}
\def\ddx{\ddot x}
\def\qu{\quad}
\def\al{\alpha}
\def\sul{\displaystyle\sum\limits}
\def\Lr{\Leftrightarrow}
\def\we{\wedge}
\def\I{\^I}
\def\a{\^a}
\def\lri{\longrightarrow}
\def\Lri{\Longrightarrow}
\def\Llr{\Longleftrightarrow}
\def\llr{\longleftrightarrow}
\def\br{R}
\def\om{\omega}
\def\qu{\quad}
\def\ld{\ldots}
\def\fo{\forall}
\def\di{\displaystyle}
\def\al{\alpha}
\def\ep{\epsilon}
\def\su{\subset}
\def\gaa{\Gamma}
\def\ty{\infty}
\def\va{\varphi}
\def\noa{\noalign}
\def\ti{\times}
\def\pp{\prime}

\setcounter{page}{1}
\textheight=19cm \textwidth=13cm
\oddsidemargin=16mm \evensidemargin=16mm
\date{}

\begin{abstract}
We illustrate the flow or wave character of the metrics and
curvatures of evolving manifolds, introducing the Riemann flow
and the Riemann wave via the bialternate product Riemannian metric. 
This kind of evolutions are new and very natural to understand certain flow or wave phenomena in the nature
as well as the geometry of evolving manifolds. It possesses many
interesting properties from both mathematical and physical point of views.

We reveal the novel features of Riemann flow and Riemann wave,
as well as their versatility, by new original results. The main results refer to: 
(1) Riemann flow PDE, (3) connection between Ricci flow and Riemann flow,
(4) the meaning of the Riemann flows on constant
curvature manifolds, (5) the gradient Riemann solitons, (6) the
infinitesimal deformations (linearizations) of Ricci flow and of Riemann flow PDEs, (7) the existence of
Ricci or Riemann curvature blow-up at finite-time singularities of the flow, 
(8) the Riemann wave PDE and its meaning on constant curvature
manifolds, (9) the existence of Ricci or Riemann curvature blow-up at finite-time singularities of the wave,
(10) the general form of some essential PDE systems on Riemannian manifolds.
\end{abstract}




\section{Motivation}

The papers of Richard Hamilton [1]-[5] and some lectures at 1990
Summer Research Institute, Differential Geometry, University of
California, Los Angeles, July 8-28, stimulated me to think about
the Riemann flow and to discuss with other mathematicians about
the subject. After years, I have introduced and studied these ideas by an open problem
in the book [18] and some explicit results in the papers [6], [21], [24]. 
Since the Hamilton Ricci flow, defined by
$$\di\frac{\pa g}{\pa t}(t,x) = - 2 \,\,\hbox{Ric}(g(t,x)),$$
was used recently by Grisha Perelman [9]-[11] to prove two very deep
theorems in topology, namely the geometrization and Poincar\' e
conjectures, I decided to bring into attention my original ideas, writing again about the
Riemann flow and the Riemann wave via the bialternate product Riemannian metric. 
The Riemann wave is connected to the Ricci wave introduced and studied recently by De-Xing
Kong and Kefeng Liu [15], [16]. 

Within the last years, I was impressed by several
researchers which have obtained original results starting from the
theory of Hamilton and Perelman (see also [7], [8], [12], [13], [14], [17]). 
Also, the Research Blog [25], which rises the problem of the curvature flow,
gave me an impulse to orient my research in this direction.  

Section 1 is a motivation for our point of view. Section 2 underlines the relation between the Riemannian metric
and the bialternate product Riemannian metric. Section 3 defines and studies the Riemann flow. The original
results include the relations between Ricci flow and Riemann flow, some existence and 
uniqueness theorems, the meaning of the Riemann flows on constant
curvature manifolds and the gradient Riemann solitons. Section 4 analyzes the
infinitesimal deformations (linearizations) of Ricci flow and of Riemann flow PDEs. Section 5 shows the existence of
Ricci or Riemann curvature blow-up at finite-time singularities of the flow.
Section 6 defines and studies the Riemann waves and their meaning on constant curvature
manifolds. Section 7 proves the existence of Ricci or Riemann curvature blow-up at finite-time singularities of the wave.
Section 8 points out the general form of some essential PDE systems on Riemannian manifolds.

\section{Bialternate product Riemannian metric}

Let $(M,g)$ be a {\it Riemannian manifold} of dimension $n$. 
The number of distinct components of the metric $g$ is $\frac{n(n + 1)}{2}$.
The Riemannian metric $g = (g_{ij})$ determines the {\it Christoffel symbols}
$$\Gamma^i_{jk} = \di\frac{1}{2}g^{il}\left(\di\frac{\pa g_{lj}}{\pa x^k} + \di\frac{\pa g_{lk}}{\pa x^j} - \di\frac{\pa g_{jk}}{\pa x^l}\right)$$
$$ =  \di\frac{1}{2}g^{im}\left(\de^r_m\, \de^s_j\, \de^t_k + \de^r_m\, \de^s_k \,\de^t_j - \de^t_m \,\de^r_j \,\de^s_k \right) \di\frac{\pa g_{rs}}{\pa x^t}$$
and this connection determines the {\it curvature (Riemannian) tensor field}
$\hbox{Riem}(g)$ of components
$$
R_{ijkl} = \di\frac{1}{2} \left(\di\frac{\pa^2 g_{ik}}{\pa x^j\pa x^l} + \di\frac{\pa^2 g_{jl}}{\pa x^i\pa x^k} -
\di\frac{\pa^2g_{jk}}{\pa x^i\pa x^l} - \di\frac{\pa^2 g_{il}}{\pa
x^j\pa x^k}\right) - g_{mn}\left( \Gamma^m_{jk}\Gamma^n_{il} -
\Gamma^m_{jl}\Gamma^n_{ik}\right)
$$
$$
= 2\, \de^p_{[i}\,\de^q_{j]}\, \de^r_{[k}\,\de^s_{l]}\,\di\frac{\pa^2 g_{pr}}{\pa x^q\pa x^s} +
2\,\de^q_j\, \de^p_i\, \de^r_{[k}\,\de^s_{l]}\,g_{mn}\Gamma^m_{qr}\,\Gamma^n_{ps}.
$$
It follows the {\it Ricci tensor field} $\hbox{Ric}(g)$ of components $R_{jk} =
g^{il}R_{ijkl}$ and the {\it curvature scalar} $R = g^{jk}R_{jk}$. The number of distinct components of the Riemannian
curvature tensor is $\frac{n^2(n^2 - 1)}{12}$.

The Riemannian metric $g$ induces the {\it bialternate product Riemannian metric}
$$G = g\odot g, \,G_{ijkl} = g_{ik}g_{jl} - g_{il}g_{jk} = 2\,\de^q_j\, \de^p_i\, \de^r_{[k}\,\de^s_{l]} \,g_{pr}\,g_{qs}$$
on $2$-forms. Then $(\Lambda^2(M), G = g\odot g)$ is the Riemannian manifold of skew-symmetric $2$-forms.
The number of distinct components of the metric $G$ is $\frac{n^2(n^2 - 1)}{12}$.
If $n=1$, then $G=0$. If $n=2$, then $$G_{1212}= - G_{1221} =
-G_{2112} = G_{2121} = \det\,(g_{ij})={\cal G}.$$ If $n\geq 3$,
then the bialternate product Riemannian metric $G$ determines the
Riemannian metric $g$ since $\hbox{rank}\,g = n\geq 3$ and
$$2g_{ij}(g_{ks}G_{mlnr} + g_{kn}G_{lmsr} + g_{kr}G_{mlsn}) = G_{mijn}G_{klrs} + G_{mijs}G_{klnr} + G_{mijr}G_{klsn}$$
$$- G_{lijn}G_{kmrs} - G_{lijs}G_{kmnr} - G_{lijr}G_{kmsn} - G_{mljs}G_{kirn} - G_{mljr}G_{kins} - G_{mljn}G_{kisr}$$
is an identity.

\section{Riemann flow via the bialternate \\product Riemannian metric}

As it was proposed in [18], let us study the flow induced by the Riemannian tensor field
$\hbox{Riem}(g(t,x))$, via bialternate product Riemannian metric
$G(t,x) = g(t,x)\odot g(t,x),\,t\in R,\,\,x\in M$, improving and developing the ideas in [6], [21], [24].

Let $M$ be a smooth closed (compact and without boundary) manifold.
Can we equip the manifold $M$ with a family of smooth Riemannian metrics $g(t,x)$
satisfying  the {\it backward heat ultrahyperbolic PDEs system}
$$\di\frac{\pa G}{\pa t}(t,x) = - 2 \,\,\hbox{Riem}(g(t,x)),\eqno (1)$$
where $(\Lambda^2(M), G(t,x) = g(t,x)\odot g(t,x))$ is the
associated Riemannian manifold of skew-symmetric $2$-forms? The {\it
Riemann flow} via the bialternate product Riemannian metric $G$ is a
mean of processing the Riemannian metric $g(t,x)$ by allowing it to
evolve under the PDE (1). This flow is driven by the Riemann tensor field. 
If $n=2$, then the PDE (1) is sub-determined; if $n =3$, then the system (1) is determined
(the number of equations is equal to the number of unknowns); if $n > 3$
then the system (1) is over-determined.

The classical (backward or forward) heat
ultrahyperbolic PDEs with standard Goursat conditions, Cauchy
conditions, Picard conditions, mixed conditions are well studied
with the help of the energy inequality method and Riemann functions.

Let ${\cal G} = \hbox{det}\,(g_{ij})$ be the determinant of the
metric $g$. Then
$$\di\frac{\pa {\cal G}}{\pa t} =\di\frac{\pa {\cal G}}{\pa g_{ik}}\di\frac{\pa g_{ik}}{\pa t}=
{\cal G}g^{ik}\di\frac{\pa g_{ik}}{\pa t},\,\,\di\frac{\pa
\hbox{ln}\,{\cal G}}{\pa t} = g^{ik}\di\frac{\pa g_{ik}}{\pa t}.$$

{\bf Theorem 3.1}. {\it If $n \geq 3$, then the Riemann type flow (1)
induces a Ricci flow of the form
$$\di\frac{\pa g_{jk}}{\pa t} = -\di\frac{2}{2-n} R_{jk}+
\di\frac{1}{2-n}\di\frac{\pa {\cal G}}{\pa t}g_{jk},\,\,\,(1 - n)\di\frac{\pa \hbox{ln}\,{\cal G}}{\pa t}=- R.$$}

{\bf Theorem 3.2}. {\it If $n \geq 3$, then the Riemann type flow
$$\di\frac{\pa G}{\pa t}(t) = \al\,\, \hbox{Riem}(g(t))+ \be\,\, \di\frac{\pa {\cal G}}{\pa t}(t)G(t),$$
with $\al = 2(n-2),\,\be  = \di\frac{1}{n-1}$, determines a standard Ricci flow
$$\di\frac{\pa g_{jk}}{\pa t} = - 2\,\, R_{jk}.$$}
In simple situations, the Ricci flow was used to deform the metric
$g(0,x)$ into a metric $g(t,x)$ distinguished by its curvature [1]-[5], [7]-[14], [17].

Conversely, we have the following statements.

{\bf Theorem 3.3}. {\it 1) On 3-dimensional Riemannian manifolds,
the Ricci flow
$$\frac {\pa g_{ij}}{\pa t} = \frac{R}{2} \,\,g_{ij} - 2\, R_{ij}$$
determines the Riemann flow (1).

2) On conformally flat Riemannian manifolds with dimension $n\geq
4$, the Ricci flow
$$\frac {\pa g_{ij}}{\pa t} =\frac{R}{2(n-1)(n-2)}\,\, g_{ij} - \frac{2}{n-2}\,\,R_{ij}$$
determines the Riemann flow (1).

3) On a Riemannian manifold with dimension $n\geq 4$, the Ricci flow
$$\frac {\pa g_{ij}}{\pa t} =\frac{R}{2(n-1)(n-2)}\,\, g_{ij} - \frac{2}{n-2}\,\,R_{ij}$$
determines a Riemann type flow
$$\frac {\pa G_{ijkl}}{\pa t} = - 2\,(R_{ijkl} - C_{ijkl}),$$
where $C_{ijkl}$ is the Weyl conformal tensor field.}

{\bf Hint}. The PDEs system
$$\frac {\pa g_{ij}}{\pa t} = \al\,\, g_{ij} + \be\,\, R_{ij}$$
implies the PDEs system
$$\frac {\pa G_{ijkl}}{\pa t} = 2\al\,\, G_{ijkl} + \be\,\, (R_{ik}g_{jl}+g_{ik}R_{jl} - R_{il}g_{jk} - g_{il}R_{jk}).$$
On the other hand,

1) in dimension three, the curvature tensor field has the local expression
$$R_{ijkl}= g_{ik}R_{jl} + g_{jl}R_{ik} - g_{il}R_{jk} - g_{jk}R_{il} -\frac{R}{2}\,\,G_{ijkl};$$

2) the curvature tensor field of a conformally flat Riemannian manifold, with dimension $n\geq 4$, is given by
$$R_{ijkl}= \frac{1}{n-2}\left(g_{ik}R_{jl} + g_{jl}R_{ik} - g_{il}R_{jk} - g_{jk}R_{il}\right) - \frac{R}{(n-1)(n-2)}\,\,G_{ijkl}.$$

3) The Weyl conformal tensor field is
$$C_{ijkl} = R_{ijkl} - \frac{1}{n-2}\left(g_{ik}R_{jl} + g_{jl}R_{ik} - g_{il}R_{jk} - g_{jk}R_{il}\right) + \frac{R}{(n-1)(n-2)}\,\,G_{ijkl}.$$

{\bf Theorem 3.4} (short time existence). {\it Given a smooth
metric $g_0$ on a closed manifold $M$, there exist $\epsilon > 0$
and a smooth family of metrics $g(t)$ for $t \in [0, \epsilon]$
such that
$$G(t) = g(t)\circ g(t), \,\di\frac{\pa G}{\pa t} = - 2\,\, \hbox{Riem}(g(t)),\,t \in [0, \epsilon],\,g(0) = g_0.$$}

{\bf Proof}. Let $\bar g_t$ be a solution to the Ricci flow $\di\frac{\pa g_{jk}}{\pa t} = - 2\,\, R_{jk}.$ Let us find
a time-dependent diffeomorphism $\varphi_t$ from our manifold $M$ to itself, with $\varphi_0 = id$, so that the metric
$g_t  = \varphi^*_t \bar g_t$ is a solution to the Riemann flow. Since
$$\frac{\pa g_t}{\pa t} = \varphi^*_t \frac{\pa \bar g_t}{\pa t} + \varphi^*_t {\cal L}_{\frac{\pa \varphi_t}{\pa t}}\bar g_t,$$
we obtain
$$
\frac{\pa G_t}{\pa t} = \varphi^*_t \frac{\pa \bar G_t}{\pa t} + \varphi^*_t {\cal L}_{\frac{\pa \varphi_t}{\pa t}}\bar G_t.
$$
On the other hand,
$$
\varphi^*_t \frac{\pa \bar G_t}{\pa t} + \varphi^*_t {\cal L}_{\frac{\pa \varphi}{\pa t}}\bar G_t = - 2\,\, \varphi^*_t (Riem( \bar g_t) - C(\bar g_t)).
$$
We select the diffeomorphism $\varphi_t$ by the condition ${\cal L}_{\frac{\pa \varphi_t}{\pa t}}\bar G_t = 2\,\,\varphi^*_t C(\bar g_t)$.

{\bf Theorem 3.5} (uniqueness of solutions). {\it Suppose
$g_1(t)$ and $g_2(t)$ are two Riemann flows on a closed manifold
$M$, for $t \in [0, \epsilon],\,\, \epsilon > 0$. If $g_1(s) =
g_2(s)$ for some $s \in [0, \epsilon]$, then $g_1(t) = g_2(t)$ for
all $t \in [0, \epsilon]$.}

Combining these two Theorems, we can talk about the Riemann
flow with initial metric $g_0$, on a maximal time interval [0, T).
In this situation, maximal means that either $T = \infty$, or that $T <
\infty$ but there do not exist $\epsilon > 0$ and a smooth Riemann flow
$\hat g(t)$ for $t \in [0, T +\epsilon)$ such that $\hat g(t) =
g(t)$ for $t \in [0, T)$.

If the initial metric $g(0,x)$ is Riemann flat, i.e.,
$\hbox{Riem}(g(0,x))=0$, then $g(t,x) = g(0,x)$ is obviously a
solution of the evolution PDE (1). Consequently, each Riemann flat
metric $g(x)$ is a steady solution of the heat ultrahyperbolic PDEs system.

\subsection{Riemann flows on constant curvature manifolds}

It is worth pointing out here that the Riemann tensor
is homogeneous under uniform scalings of the metric. More precisely, if $g(t) = f(t)g_0$, $f(t) > 0$,
$f(0) = 1$, then $\hbox{Riem}(g(t)) = f(t)\,\hbox{Riem}(g_0)$.

If we take a metric $g_0$ such that $\hbox{Riem}(g_0) = \la\,\, G_0$
for some constant $\la \in R$ (these metrics are known as constant
curvature metrics), then a solution $g(t)$ of
PDE (1) with $g(0) = g_0$ is given by $g(t) = (1 - \la t)\,g_0$ and hence
$G(t) = (1 - \la t)^2\,G_0$. In particular, for the round unit sphere $(S^n, g_0)$, $n\geq 2$, we have
$\hbox{Riem}(g_0) = (n-1)G_0$, so the evolution metric is $g(t) =
(1- (n-1)t)\,\,g_0$ and the sphere collapses to a point at
finite time $T = \di\frac{1}{n - 1}$. An alternative example of
this type would be if $g_0$ were a hyperbolic metric, that is, of
constant sectional curvature $-1$. In this case, for $n\geq 2$, we have
$\hbox{Riem}(g_0) = - (n -1)G_0$. Consequently the evolution metric is $g(t) =
(1 + (n - 1)t)\,g_0$ and the manifold expands homothetically for all time. 

Of course a metric is flat, i.e., the Riemannian
manifold is locally isometric to the Euclidean space, if and only
if the associated Riemann curvature tensor vanishes. These metrics
can be regarded as the fixed points of the Riemann flow.

\subsection{Gradient Riemann soliton}

We can regard as {\it generalized fixed points} of the Riemann flow
those manifolds which change only by a diffeomorphism and a
rescaling under the Riemann flow. To simplify, instead of $g(t,x)$
we shall write $g(t)$. Let $(M,g(t))$ be a solution of the Riemann
flow, and suppose $\varphi_t: M\to M,\,\varphi_0 = id$ is a
time-dependent family of diffeomorphisms and $\si (t)$ is a
time-dependent scale factor (with $\si (0)=1$). If we have $g(t)
=\si (t)\varphi_t^*\,g(0)$, then the solution $(M, g(t))$ or
$(\Lambda^2(M), G(t) = g(t)\odot g(t),\,t\in R)$ is called {\it
Riemann soliton}. Of course
$$G_{ijkl} (t)= g_{ik}(t)g_{jl}(t) - g_{il}(t)g_{jk}(t)$$
$$ = \si^2 (t)\varphi_t^*\,g_{ik}(0)\varphi_t^*\,g_{jl}(0) - \si^2 (t)\varphi_t^*\,g_{il}(0)\varphi_t^*\,g_{jk}(0).$$
Computing the derivative $\di\frac{\pa G}{\pa t}(t)$ and
evaluating it at $t=0$, we find
$$\di\frac{\pa G}{\pa t}(t)|_{t=0} = 2\si(t) \di\frac{d\si}{dt}(t)\varphi_t^*g(t)\circ \varphi_t^*g(t)|_{t=0}$$
$$+ \si^2(t)\left(\di\frac{\pa}{\pa t}\varphi_t^*g(t)\circ \varphi_t^*g(t) + \varphi_t^*g(t)\circ \di\frac{\pa}{\pa t}\varphi_t^*g(t) \right)|_{t=0}$$
$$ =2 \si^{\prime}(0)\,g(0)\circ g(0) +{\cal L}_V g(0)\circ g(0) + g(0)\circ {\cal L}_V g(0)= -2\, \hbox{Riem}(g(0)),$$
where $V = \di\frac{\pa\varphi_t}{\pa t}$, by the definition of
{\it Lie derivative}. We replace $\si^{\prime}(0) = \la$ and we
write in coordinates
$$-2 R_{ijkl} = 2\la (g_{ik}(0)g_{jl}(0) - g_{il}(0)g_{jk}(0))$$
$$+ {\cal L}_Vg_{ik}(0)g_{jl}(0)+ g_{ik}(0){\cal L}_Vg_{jl}(0) - {\cal L}_Vg_{il}(0)g_{jk}(0) - g_{il}(0){\cal L}_Vg_{jk}(0)$$
$$= 2\la (g_{ik}(0)g_{jl}(0) - g_{il}(0)g_{jk}(0))$$
$$+ (\nabla_i V_k +\nabla_k V_i)g_{jl}(0)+ g_{ik}(0)(\nabla_j V_l +\nabla_l V_j)$$
$$ - (\nabla_i V_l +\nabla_l V_i)g_{jk}(0) - g_{il}(0)(\nabla_j V_k +\nabla_k V_j).$$
Suppose $V = \nabla f$. Then the foregoing relation becomes
$$R_{ijkl} + \la (g_{ik}(0)g_{jl}(0) - g_{il}(0)g_{jk}(0))$$
$$+ g_{jl}(0) \nabla_i \nabla_k f + g_{ik}(0)\nabla_j \nabla_l f - g_{jk}(0) \nabla_i \nabla_l f - g_{il}(0)\nabla_j \nabla_k f =0.$$
Such solutions we call {\it gradient Riemann solitons}. A gradient
Riemann soliton is called {\it shrinking} if $\la <0$, {\it
static} if $\la =0$, and {\it expanding} if $\la > 0$.

\section{Infinitesimal deformations}

\subsection{Linearizing Ricci flow PDE}

Let $g(t,x),\,\,t\in I\subset R$, be a Ricci flow. Let $g(t,x;\epsilon)$ be a differentiable variation of $g(t,x)$, i.e.,
$$\di\frac{\pa g}{\pa t}(t,x;\epsilon) = - 2\,\, \hbox{Ric}(g(t,x;\epsilon)),\, g(t,x;0) = g(t,x).$$
We take the partial derivative with respect to $\epsilon$ and we denote $\di\frac{\pa g}{\pa \epsilon}|_{\epsilon = 0} = h$.
The relations $g_{ij}g^{jk} = \de^k_i$ implies
$$\frac{\pa g^{jk}}{\pa \epsilon}|_{\epsilon = 0} = - h_{ij}g^{jk}g^{il} = - h^{lk}.$$
We find the {\it single-time infinitesimal deformation PDE}
$$
\frac{\pa h}{\pa t} = - 2\left( \frac{\pa \hbox{Ric}}{\pa g^{ij}} (- h^{ij})+\frac{\pa \hbox{Ric}}{\pa {\frac{\pa g_{kl}}{\pa x^i}}} \frac{\pa h_{kl}}{\pa x^i}
 + \frac{\pa \hbox{Ric}}{\pa {\frac{\pa^2 g_{kl}}{\pa x^i \pa x^j}}} \frac{\pa^2 h_{kl}}{\pa x^i\pa x^j}\right)
$$
around a solution $g(t,x)$ of the Ricci PDE. Of course, this new PDE is the linearization of the Ricci PDE along a solution $g(t,x)$.

\subsection{Linearizing Riemann flow PDE}

Let $g(t,x),\,\,t\in I\subset R$, be a Riemann flow. Let $g(t,x;\epsilon)$ be a differentiable variation of $g(t,x)$, i.e.,
$$\di\frac{\pa G}{\pa t}(t,x;\epsilon) = - 2\,\, \hbox{Riem}(g(t,x;\epsilon)),\, g(t,x;0) = g(t,x).$$
We take the partial derivative with respect to $\epsilon$ and we denote $\di\frac{\pa g}{\pa \epsilon}|_{\epsilon = 0} = h$.
The relations $g_{ij}g^{jk} = \de^k_i$ implies
$$\frac{\pa g^{jk}}{\pa \epsilon}|_{\epsilon = 0} = - h_{ij}g^{jk}g^{il} = - h^{lk}.$$
We find the {\it single-time infinitesimal deformation PDE}
$$
\frac{\pa }{\pa t} (h\odot g + g\odot h) = - 2\left( \frac{\pa \hbox{Riem}}{\pa g^{ij}}
(- h^{ij})+\frac{\pa \hbox{Riem}}{\pa {\frac{\pa g_{kl}}{\pa x^i}}} \frac{\pa h_{kl}}{\pa x^i}
 + \frac{\pa \hbox{Riem}}{\pa {\frac{\pa^2 g_{kl}}{\pa x^i \pa x^j}}} \frac{\pa^2 h_{kl}}{\pa x^i\pa x^j}\right)
$$
around a solution $g(x,t)$ of the Riemann PDE. Of course, this new PDE is the linearization of the Riemann PDE along a solution $g(t,x)$.

\section{Blow-ups at finite-time singularities for \\geometric flows}

\subsection{Ricci estimation and Ricci curvature blow-up at \\finite-time singularities}

This subsection will be devoted to proving that, if the Ricci flow $g(t,x)$, defined by
$$\di\frac{\pa g}{\pa t}(t,x) = - 2 \,\,\hbox{Ric}(g(t,x)),$$
becomes singular in finite time, then the Ricci curvature must explode as we approach the singular time $T$.

{\bf Lemma 4.1} (Riemannian metric equivalence). {\it If $g(t)$ is a Ricci flow for $t \in [0, T)$ and
$||\hbox{Ric}|| < m$ on $M\times [0, T)$, then}
$$e^{-2mt} g(0) \leq g(t) \leq e^{2mt} g(0), \forall t\in [0,T),$$
$$||g(t)||^2 - ||g(0)||^2 \leq C\, e^{4m^2t}, \forall t\in [0,T).$$

{\bf Remark}. Given two symmetric $(0, 2)$-tensors $a$ and $b$,
we write $a \geq b$ to mean that $a - b$ is positive semidefinite.

{\bf Proof}. Since $\frac{\pa}{\pa t}\,g(X,X) = - 2\, \hbox{Ric}(X,X)$, we have
$$\left|\frac{\pa}{\pa t}\,g(X,X)\right| = 2 ||\hbox{Ric}|| g(X,X)$$
and hence
$$\left|\frac{\pa}{\pa t}\,\ln g(X,X)\right| \leq 2 ||\hbox{Ric}||,$$
for any nonzero tangent vector $X$. In this way,
$$\left|\ln \frac{g(t)(X,X)}{g(0)(X,X)}\right| \leq 2 mt.$$

To obtain the second formula, we remark that
$$\di\frac{1}{2}\frac{\pa}{\pa t} ||g||^2 = - 2 <\hbox{Ric}, g> \,\,\leq 2\,\,||\hbox{Ric}||^2\,||g||^2$$
implies
$$||g(t)||^2 - ||g(0)||^2 \leq C\, e^{4\di\int_0^t ||\hbox{Ric}||^2 ds}.$$

{\bf Theorem 4.2} {\it Let $M$ be a closed manifold. If $g(t)$ is a Ricci flow on a maximal
time interval $[0, T)$ and $T < \infty$, then}
$$\lim_{t\nearrow T}\left(\sup_{x \in M} ||\hbox{Ric}(t,x)||\right) = \infty.$$

{\bf Proof}. We will deal with the contrapositive of this theorem.
Let us show that if $||\hbox{Ric}(t,x)||\leq m,\, t\in [0,T)$, then the Ricci flow $g(t)$
can be extended to a larger time interval $[0, T + \ep)$.
Indeed, the relation
$$g(t) = g(0) - 2\int_0^t\, \hbox{Ric}(s,x)ds,\,\, t\in [0,T)$$
implies
$$||g(t^{\prime}) - g(t^{\prime\prime})|| \leq 2m |t^{\prime} - t^{\prime\prime}|,\,\, \forall t^{\prime}, t^{\prime\prime}\in [0,T).$$
The Cauchy Criterion shows that $\lim_{t\nearrow T}\,g(t)$ exists, while $\lim_{t\nearrow T}\,\frac{\pa g}{\pa t}(t)$
exists since $\lim_{t\nearrow T}\,\hbox{Ric}(t,x)$ exists
(see the definition of the Ricci flow). Consequently $g(t),\,t \in [0,T]$ is a solution for the Ricci flow.

In this case, the Ricci flow $g(t)$ may be extended from being a smooth
solution on $[0, T)$ to a smooth solution on $[0, T]$. Then we take $g(T)$ to
be an initial metric in our short-time existence theorem (Theorem 3.4)
in order to extend the flow to a Ricci flow for $t \in [0, T + \ep)$, contradicting
the assumption that $[0, T)$ is a maximal time interval.

\subsection{Riemann estimation and curvature blow-up at \\finite-time singularities}

This subsection will be devoted to proving that, if the Riemann flow becomes singular in finite time, then the
curvature must explode as we approach the singular time $T$.

{\bf Lemma 4.3} (bialternate product Riemannian metric equivalence). {\it If $g(t)$ is a Riemann flow for $t \in [0, T)$ and
$||\hbox{Riem}|| < m$ on $M\times [0, T)$, then}
$$e^{-2mt} G(0) \leq G(t) \leq e^{2mt} G(0), \forall t\in [0,T),$$
$$||G(t)||^2 - ||G(0)||^2 \leq C\, e^{4m^2t}, \forall t\in [0,T).$$

{\bf Proof}. Since $\frac{\pa}{\pa t}\,G(X,Y,X,Y) = - 2\, \hbox{Riem}(X,Y,X,Y)$, we have
$$\left|\frac{\pa}{\pa t}\,G(X,Y,X,Y)\right| = 2\, ||\hbox{Riem}|| G(X,Y,X,Y)$$
and hence
$$\left|\frac{\pa}{\pa t}\,\ln G(X,Y,X,Y)\right| \leq 2\, ||\hbox{Riem}||,$$
for any nonzero tangent vectors $X, Y$. In this way,
$$\left|\ln \frac{G(t)(X,Y,X,Y)}{G(0)(X,Y,X,Y)}\right| \leq 2 mt.$$

To obtain the second formula, we remark that
$$\di\frac{1}{2}\frac{\pa}{\pa t} ||G||^2 = - 2 <\hbox{Riem}, G> \,\,\leq 2\,\,||\hbox{Riem}||^2\,||G||^2$$
implies
$$||G(t)||^2 - ||G(0)||^2 \leq C\, e^{4\int_0^t ||\hbox{Riem}||^2 ds}.$$

{\bf Theorem 4.4} {\it Let $M$ be a closed manifold with the dimension $n \geq 3$. If $g(t)$
is a Riemann flow on a maximal time interval $[0, T)$ and $T < \infty$, then}
$$\lim_{t\nearrow T}\left(\sup_{x \in M} ||\hbox{Riem}(t,x)||\right) = \infty.$$

{\bf Proof}. Let us show that if $||\hbox{Riem}(t,x)||\leq m,\, t\in [0,T)$, then the Riemann
flow $g(t)$ can be extended to a larger time interval $[0, T + \ep)$. Indeed, the relation
$$G(t) = G(0) - 2\int_0^t\, \hbox{Riem}(s,x)ds,\,\, t\in [0,T)$$
implies
$$||G(t^{\prime}) - G(t^{\prime\prime})|| \leq 2m |t^{\prime} - t^{\prime\prime}|,\,\, \forall t^{\prime}, t^{\prime\prime}\in [0,T).$$
The Cauchy Criterion shows that $\lim_{t\nearrow T}\,G(t)$ exists, while $\lim_{t\nearrow T}\,\frac{\pa G}{\pa t}(t)$
exists since $\lim_{t\nearrow T}\,\hbox{Riem}(t,x)$ exists
(see the definition of the Riemann flow). Consequently $g(t),\,t \in [0,T]$ is a solution for the Riemann flow.

In this case, the Riemann flow $g(t)$ may be extended from being a smooth
solution on $[0, T)$ to a smooth solution on $[0, T]$. Then we take $g(T)$ to
be an initial metric in our short-time existence theorem (Theorem 3.4)
in order to extend the flow to a Riemann flow for $t \in [0, T + \ep)$, contradicting
the assumption that $[0, T)$ is a maximal time interval.

\section{Riemann wave via \\the bialternate product Riemannian metric}

Let us study the wave of the Riemannian tensor field
$\hbox{Riem}(g(t,x))$, via bialternate product Riemannian metric
$G(t,x) = g(t,x)\odot g(t,x),\,t\in R,\,\,x\in M$, improving and developing the ideas in  [6], [21], [24].

Let $M$ be a smooth closed (compact and without boundary) manifold.
Can we equip the manifold $M$ with a family of smooth Riemannian metrics $g(t,x)$
satisfying  the {\it wave ultrahyperbolic PDEs system}
$$\di\frac{\pa^2 G}{\pa t^2}(t,x) = - 2 \,\,\hbox{Riem}(g(t,x)),\eqno (2)$$
where $(\Lambda^2(M), G(t,x) = g(t,x)\odot g(t,x))$ is the associated
Riemannian manifold of skew-symmetric $2$-forms?
The {\it Riemann wave} via the bialternate product Riemannian metric $G$ is a mean of processing the Riemannian metric $g(t,x)$ by
allowing it to evolve under the PDE (2). If $n=2$, then the PDE (2) is sub-determined; if $n =3$, then the system (1) is determined
(the number of equations is equal to the number of unknowns); if $n > 3$
then the system (1) is over-determined.

{\bf Theorem 4.1} (uniqueness of solution). {\it Let $(M,
g_0(x))$ be a compact Riemannian manifold and $k_1(x)$ be a
$(0,2)$ symmetric tensor field on $M$. Then there exists a
constant $\epsilon > 0$ such that the initial value problem
$$\di\frac{\pa^2 G}{\pa t^2}(t,x) = -2\,\, \hbox{Riem}(g(t,x)),\,\,g(0,x)=g_0(x),\,\,\frac{\pa g}{\pa t}(0,x) = k_1(x)$$
has a unique smooth solution $g(t,x)$ on $M\times[0,\epsilon]$.}

{\bf Theorem 4.2}. {\it If $(M,\,\,
g_{ij}(0,x)=u_0(x^1)\delta_{ij},\,u_0(x^1)>0)$ is a conformally
flat Riemannian manifold, then the Cauchy problem
$$\di\frac{\pa^2 G}{\pa t^2}(t,x) = - 2\,\, \hbox{Riem}(g(t,x)),\,\,g(0,x)=g_0(x),\,\,\frac{\pa g}{\pa t}(0,x) = k_1(x)$$
has a unique smooth solution for all time $t\in R$, and the
solution metric has the form
$g_{ij}=u(t,x^1)\delta_{ij},\,u(t,x^1)>0$.}

{\bf Hint}. Let (M, g) be a conformally flat Riemannian manifold.
If we take a Riemannian metric with the components $g_{ii} =
H_i^2,\,g_{ij}=0,\,i\not= j$ (see Lam\' e coefficients), then
$$R_{hijk}=0,\,\, (h,i,j,k\not=)$$
$$R_{hiik}= - H_i\left(\frac{\pa^2 H_i}{\pa x^h\pa x^k} - \frac{\pa H_i}{\pa x^h}\frac{\pa \ln \,H_h}{\pa x^k} -
\frac{\pa H_i}{\pa x^k}\frac{\pa \ln \,H_k}{\pa x^h}\right),\,\,(h,i,k\not=)$$
$$R_{hiih}= - H_hH_i\left(\frac{\pa}{\pa x^h}\left(\frac{1}{H_h}\frac{\pa H_i}{\pa x^h}\right)+
\frac{\pa}{\pa x^i}\left(\frac{1}{H_i}\frac{\pa H_h}{\pa x^i}\right)+
{\sum_l{}}^{\prime}\frac{1}{H_l^2}\frac{\pa H_h}{\pa x^l}\frac{\pa
H_i}{\pa x^l}\right),$$ where the index of summation $l$ indicates
the sum for the values $1,...,n$ excepting $h$ and $i$. Replacing
$g_{ij}=u(t,x^1)\delta_{ij}$, we find
$$\left(\frac{\pa u}{\pa t}\right)^2 +\left(\frac{\pa u}{\pa x^1}\right)^2 + u\left(\frac{\pa^2 u}{\pa t^2}- \frac{\pa^2 u}{\pa {x^1}^2}\right)=0,$$
with the initial conditions $u(0,x^1)=u_0(x^1),\,\,\di\frac{\pa
u}{\pa t}(0,x^1) = u_1(x^1).$

\subsection{Riemann waves on constant curvature manifolds}

We take a metric $g_0$ such that $\hbox{Riem}(g_0) = \la G_0$
for some constant $\la \in R$ (these metrics are known as constant
curvature metrics). Then a solution $g(t)$ of
PDE (2) with $g(0) = g_0$ is of the form $g(t) = f(t)\,g_0,\,f(t) >0, f(0) = 1, f^{\prime}(0) = v$,
and hence $G(t)$ is of the form $G(t) = f^2(t)G_0$, if and only if
$${f^{\prime}}^2(t) + f(t) f^{\prime\prime}(t) + \la f(t) = 0,\, f(0) = 1, f^{\prime}(0) = v.$$

If $\la < 0$, the solution is the polynomial function
$$f(t) = 1 + v t  - \frac{\la}{6}\, t^2.$$
In particular, let $g_0$ be a hyperbolic metric, that is, of constant sectional curvature
$-1$. In this case, $n\geq 2$, $\hbox{Riem}(g_0) = -(n -1)G_0$,
the evolution metric is $g(t) = (1 + v t  - \frac{\la}{6}\, t^2)\,g_0$ and the manifold
expands homothetically for all time. 

If $\la > 0$, there exists $T$ depending on $v$ such that the solution 
$f(t),\, t \in [0,T)$ is a concave function with $\lim_{t\to T} f(t) = 0$.
In particular, for the round unit sphere $(S^n, g_0)$, $n\geq 2$, we have
$\hbox{Riem}(g_0) = (n-1)G_0$, so the evolution metric is $g(t) =
f(t)\,g_0$ and the sphere collapses to a point when $t\to T$. 

If the initial metric $g(0,x)$ is Riemann flat, i.e., $\hbox{Riem}(g(0,x))=0$, then
$g(t,x) = g(0,x)$ is obviously a solution of the evolution PDE
(2). Consequently, each Riemann flat metric $g(x)$ is a steady
solution of the wave ultrahyperbolic PDEs system.

{\bf Open problem}. The case of multitime Ricci wave or multitime
Riemann wave PDEs does not seem to have been investigated. 
The study of this case is caused not only by theoretical interest, 
but also by practical necessity of describing the deformations (see
[18]-[23] for inspiration).

\section{Blow-ups at finite-time singularities for \\geometric waves}

Let us prove that, if the Ricci wave becomes singular in finite time, then the
Ricci curvature must explode as we approach the singular time $T$.

{\bf Theorem 6.1} {\it Let $M$ be a closed manifold. If $g(t)$ is a Ricci wave, i.e.,
$$\di\frac{\pa^2 g}{\pa t^2}(t,x) = - 2 \,\,\hbox{Ric}(g(t,x)),$$
on a maximal time interval $[0, T)$ and $T < \infty$, then}
$$\lim_{t\nearrow T}\left(\sup_{x \in M} ||\hbox{Ric}(t,x)||\right) = \infty.$$

{\bf Proof}. Let us show that if $||\hbox{Ric}(t,x)||\leq m,\, t\in [0,T)$, then
the Ricci wave $g(t)$ can be extended to a larger time interval $[0, T + \ep)$.
Indeed, the relations
$$\frac{\pa g}{\pa t}(t) = \frac{\pa g}{\pa t}(0) - 2\,\int_0^t\, \hbox{Ric}(s,x)ds,\,\, t\in [0,T),$$
$$g(t) = g(0) + \frac{\pa g}{\pa t}(0)\, t - 2\, \int_0^t \int_0^u\, \hbox{Ric}(s,x)dsdu,\,\, t\in [0,T)$$
imply
$$g(t^{\prime}) - g(t^{\prime\prime}) = \frac{\pa g}{\pa t}(0)\, (t^{\prime} - t^{\prime\prime}) +
2\, \int_{t^{\prime}}^{t^{\prime\prime}} \!\!\int_0^u\, \hbox{Ric}(s,x)dsdu,\,\, t\in [0,T)$$
and hence
$$
||g(t^{\prime}) - g(t^{\prime\prime}|| \leq \left( || \frac{\pa g}{\pa t}(0)|| + 2mT\right)
|t^{\prime} - t^{\prime\prime}|,\,\, \forall t^{\prime}, t^{\prime\prime}\in [0,T).
$$
The Cauchy Criterion shows that $\lim_{t\nearrow T}\,g(t)$ exists, while $\lim_{t\nearrow T}\,\frac{\pa g}{\pa t}(t)$ and
$\lim_{t\nearrow T}\,\frac{\pa^2 g}{\pa t^2}(t)$ exist since $\lim_{t\nearrow T}\,\hbox{Ric}(t,x)$ exists
(see the definition of the Ricci wave). Consequently $g(t),\,t \in [0,T]$ is a Ricci wave.

In this case, the Ricci wave $g(t)$ may be extended from being a smooth
solution on $[0, T)$ to a smooth solution on $[0, T]$. Then we take $g(T),\,\frac{\pa g}{\pa t}(T)$ to
be an initial metric in a short-time existence theorem
in order to extend the wave to a Ricci wave for $t \in [0, T + \ep)$, contradicting
the assumption that $[0, T)$ is a maximal time interval.

Similarly, let us prove that, if the Riemann wave becomes singular in finite time, then the
Riemann curvature must explode as we approach the singular time $T$.

{\bf Theorem 6.2} {\it Let $M$ be a closed manifold with the dimension $n \geq 3$. If $g(t)$ is a Riemann wave, i.e.,
$$\di\frac{\pa^2 G}{\pa t^2}(t,x) = - 2 \,\,\hbox{Riem}(g(t,x)),$$
on a maximal time interval $[0, T)$ and $T < \infty$, then}
$$\lim_{t\nearrow T}\left(\sup_{x \in M} ||\hbox{Riem}(t,x)||\right) = \infty.$$

{\bf Proof}. Let us show that if $||\hbox{Riem}(t,x)||\leq m,\, t\in [0,T)$, then the
Riemann wave $g(t)$ can be extended to a larger time interval $[0, T + \ep)$.
Indeed, the relations
$$\frac{\pa G}{\pa t}(t) = \frac{\pa G}{\pa t}(0) - 2\,\int_0^t\, \hbox{Riem}(s,x)ds,\,\, t\in [0,T),$$
$$G(t) = G(0) + \frac{\pa G}{\pa t}(0)\, t - 2\, \int_0^t \int_0^u\, \hbox{Riem}(s,x)dsdu,\,\, t\in [0,T)$$
imply
$$G(t^{\prime}) - G(t^{\prime\prime}) = \frac{\pa g}{\pa t}(0)\, (t^{\prime} - t^{\prime\prime}) +
2\, \int_{t^{\prime}}^{t^{\prime\prime}}\!\! \int_0^u\, \hbox{Riem}(s,x)dsdu,\,\, t\in [0,T)$$
and hence
$$
||G(t^{\prime}) - G(t^{\prime\prime}|| \leq \left( || \frac{\pa G}{\pa t}(0)|| + 2mT\right)
|t^{\prime} - t^{\prime\prime}|,\,\, \forall t^{\prime}, t^{\prime\prime}\in [0,T).
$$
The Cauchy Criterion shows that $\lim_{t\nearrow T}\,G(t)$ exists, while $\lim_{t\nearrow T}\,\frac{\pa G}{\pa t}(t)$
and $\lim_{t\nearrow T}\,\frac{\pa^2 G}{\pa t^2}(t)$ exist since $\lim_{t\nearrow T}\,\hbox{Riem}(t,x)$ exists
(see the definition of the Riemann wave). Consequently $g(t),\,t \in [0,T]$ is a Riemann wave.

In this case, the Riemann wave $g(t)$ may be extended from being a smooth
solution on $[0, T)$ to a smooth solution on $[0, T]$. Then we take $g(T),\,\frac{\pa g}{\pa t}(T)$ to
be an initial metric in a short-time existence theorem
in order to extend the wave to a Riemann wave for $t \in [0, T + \ep)$, contradicting
the assumption that $[0, T)$ is a maximal time interval.

\section{Conclusions}

On Riemannian manifolds we can discuss about the ultrahyperbolic
PDEs system (constant curvature manifolds) $R_{ijkl} = \la\,\,G_{ijkl}.$
On Riemannian manifolds whose metric depends on the
evolution parameter $t$, we can add: heat ultrahyperbolic PDEs
system (Riemann flow), $\di\frac{\pa G_{ijkl}}{\pa t} = - 2\,\,
R_{ijkl};$ wave ultrahyperbolic PDEs system (Riemann wave),
$\di\frac{\pa^2 G_{ijkl}}{\pa t^2} = - 2 \,\,R_{ijkl}.$ In fact, the Riemann flow and the
Riemann wave are PDE systems for evolving the Riemannian metric to make it "rounder",
in the hope that one may draw topological conclusions from the existence of such "round" metrics.
The general form of these PDE systems is
$$\al_{ijkl}\,\,\frac{\pa^2 G_{ijkl}}{\pa t^2} + \be_{ijkl}\,\,\frac{\pa G_{ijkl}}{\pa t}+\ga_{ijkl}\,\,G_{ijkl}+\de_{ijkl}\,\, R_{ijkl}=0,$$
where the coefficients are smooth functions depending on $t,\,x$, having the same symmetries as
the curvature tensor. All these determine new classes of Riemannian metrics.

There are more ingenious ways to define such PDEs and PDIs using
spacetimes (see [6]).

\section*{Acknowledgement}
Partially supported by University Politehnica of Bucharest and by
Academy of Romanian Scientists, Bucharest, Romania.



\begin{thebibliography}{99}

\bibitem{[1]} R. Hamilton, {\it Three-manifolds with positive Ricci curvature}, J. Diff. Geo., {17} (1982), 255--306.

\bibitem{[2]} R. Hamilton, {\it Four-manifolds with positive curvature operator}, J. Diff. Geo., {24} (1986), 153--179.

\bibitem{[3]} R. Hamilton, {\it A compactness property for solutions of the Ricci flow}, American Journal of Mathematics, 117 (1992), 545--572.

\bibitem{[4]} R. Hamilton, {\it The formation of singularities in the Ricci flow}, Surveys in Differential Geometry, 2 (1995), 1--136.

\bibitem{[5]} R. Hamilton, {\it Non-singular solutions of the Ricci flow on three-manifolds}, Comm. Anal. Geom., 7 (1999), 695--729.

\bibitem{[6]} I. E. Hiric\u a, C. Udri\c ste, {\it Basic evolution PDEs in Riemannian geometry}, 
Balkan J. Geom. Appl., 17, 1 (2012), 30-40.


\bibitem{[7]} C. Hopper, B. Andrews, {\it The Ricci Flow in Riemannian Geometry}, Springer, 2009.

\bibitem{[8]} D.  M\' aximo, {\it Non-negative Ricci curvature on closed manifolds under Ricci flow},
Proceedings of the American Mathematical Society, 139, 2 (2011), 675-685.

\bibitem{[9]} G. Perelman, {\it The entropy formula for the Ricci flow and its geometric applications},
http://arXiv.org/math.DG/0211159v1, 2002.

\bibitem{[10]} G. Perelman, {\it Finite extinction time for the solutions to the Ricci flow on certain three-manifolds},
http://arXiv.org/math.DG/0307245v1, 2002.

\bibitem{[11]} G. Perelman, {\it Ricci flow with surgery on three-manifolds}, http://arXiv.org/math.DG/0303109v1, 2003.

\bibitem{[12]} N. Sheridan, {\it Hamilton's Ricci Flow}, PhD Thesis, The University of Melbourne,
Department of Mathematics and Statistics, 2006.

\bibitem{[13]} J. Porti, {\it Hamilton-Ricci flow on three dimensional manifolds},
Fall Workshop in Geometry and Physics, Bilbao, 2005.

\bibitem{[14]} B. Chow, D. Knopf, {\it The Ricci Flow: An Introduction},
AMS Math. Surveys and Monographs, 2004.

\bibitem{[15]} D.-X. Kong, K. Liu, {\it Wave character of metrics and hyperbolic geometric flow},
http://www.cms.zju.edu.cn/UploadFiles/AttachFiles/ 200682885946597.pdf

\bibitem{[16]} D.-X. Kong, K. Liu, D.-L. Xu, {\it The hyperbolic geometric flow on Riemann surfaces},
arXiv:0709.1607v3 [math.DG] 9 Jan 2008.

\bibitem{[17]} P. Topping, {\it Lectures on the Ricci flow}, 2006.

\bibitem{[18]} C. Udri\c ste, M. Ferrara, D. Opri\c s, {\it Economic Geometric Dynamics},
Monographs and Textbooks 6, Geometry Balkan Press, Bucharest, 2004.

\bibitem{[19]} C. Udri\c ste, {\it Simplified multitime maximum principle},
Balkan J. Geom. Appl., 14, 1 (2009), 102-119.

\bibitem{[20]} C. Udri\c ste, {\it Nonholonomic approach of multitime maximum principle},
Balkan J. Geom. Appl., 14, 2 (2009), 111-126. 

\bibitem{[21]} C. Udri\c ste, {\it Riemann flow and Riemann wave}, An. Univ. Vest,
Timi\c soara, Ser. Mat.-Inf., 48, 1-2 (2010), 265-274.

\bibitem{[22]} C. Udri\c ste, {\it Multitime maximum principle for curvilinear integral cost},
Balkan J. Geom. Appl., 16, 1 (2011), 128-149.

\bibitem{[23]} C. Udri\c ste, A. Bejenaru, {\it Multitime optimal control with area integral costs on boundary}, 
Balkan J. Geom. Appl., 16, 2 (2011), 138-154. 

\bibitem{[24]} C. Udri\c ste, I. \c Tevy, {\it Sturm-Liouville operator controlled by sectional curvature on Riemannian manifolds},
Balkan J. Geom. Appl., 17, 2 (2012), 129-140.

\bibitem{[25]} Research Blog 2/28/03, www.math.uic.edu/~agol/blog/030228.pdf



\end{thebibliography}
\end{document}